\documentclass[a4paper, 12pt]{amsproc}
\usepackage{mathrsfs}
\usepackage{amssymb}
\usepackage{amsmath}
\title[Steinberg modules of Reductive Groups]
{Irreducibility of Infinite Dimensional Steinberg Modules of Reductive Groups with Frobenius Maps}

\author[Ruotao Yang]{ Ruotao Yang$^{*}$}
\address{$^{*}$
Department of Mathematics\\
University of Science and Technology of China\\
Hefei\\
China} \email{yrt@mail.ustc.edu.cn}

\thanks{Supported in part by a grant of the National Natural Science Foundation of China (No. 11321101).}

\begin{document}
\baselineskip=18pt
\begin{abstract}
 Let $G$ be a connected reductive group over an algebraic closure $\bar{\mathbb{F}}_q$ of a finite field ${\mathbb{F}_q}$. In this paper it is proved that the infinite dimensional Steinberg module of $kG$ defined by N. Xi in 2014 is irreducible when $k$ is a field of positive characteristic and char$k\ne $char$\mathbb{F}_q$. For certain special linear groups, we show that the Steinberg modules of the groups are not quasi-finite with respect to some natural quasi-finite sequences of the groups.

\end{abstract}

\maketitle

\def\Cal{\mathcal}
\def\bold{\mathbf}
\def\ca{\mathcal A}
\def\cdz{\mathcal D_0}
\def\cd{\mathcal D}
\def\cdo{\mathcal D_1}
\def\bold{\mathbf}
\def\l{\lambda}
\def\le{\leq}

N. Xi studied some induced representations of infinite reductive groups with Frobenius maps (see [X]). In particular, he  defined Steinberg modules for any reductive groups by extending Steinberg's construction of Steinberg modules for finite reductive groups. These Steinberg modules are infinite dimensional when the reductive groups are infinite.

Let $G$ be a connected reductive group over the algebraic closure $\bar{\mathbb{F}}_q$ of a finite field ${\mathbb{F}_q}$ and  $k$   a  field. Xi proved that the Steinberg module of the group algebra $kG$ of $G$ over $k$ is irreducible if $k$ is the field of complex numbers or $k=\bar{\mathbb{F}}_q$ ( In fact, his proof works when char k=0 or char $\mathbb{F}_q$ ) . In this paper we prove that if $k$ has positive characteristic and char$k\ne $ char $\mathbb{F}_q$, then the Steinberg module of $kG$ remains irreducible (see Theorem 2.2).

The reductive group $G$ is quasi-finite in the sense of [X, 1.8]. For quasi-finite groups Xi  introduced the concept of quasi-finite irreducible module and raised the question whether an irreducible $kG$-module is always quasi-finite. For certain special linear groups, we show that the Steinberg modules of the groups are not quasi-finite with respect to some natural quasi-finite sequences of the groups (see Proposition 3.2).

\section{Preliminaries}
\def\ind{{\text {Ind}}}
\def\Hom{{\text {Hom}}}
\def\res{{\text {Res}}}

\def\bbf{\mathbb{F}_q}
 \def\bbfa{\mathbb{F}_{q^a}}
 \def\bbbf{\bar{\mathbb{F}}_q}
 \def\bbc{\mathbb{C}}

In this section we recall some basic facts for reductive groups defined over a finite field, for details we refer to [C].

\subsection{} Let  $G$ be  a connected reductive group over an algebraically closure $\bbbf$
 of a finite field $\mathbb{F}_q$ of $q$ elements, where $q$ is a power of a prime $p$. Assume that $G$ is defined over $\mathbb{F}_q$. Then
  $G$ has a Borel subgroup $B$ defined over $\bbf$ and $B$ contains
 a maximal torus $T$ defined over $\bbf$. The unipotent radical $U$ of $B$ is defined over $\bbf$. For any power $q^a$ of $q$,
 we denote by $G_{q^a}$ the $\mathbb{F}_{q^a}$-points of $G$ and shall identify $G$ with its $\bar{\mathbb{F}}_q$-points. Then we have $G=\bigcup_{a=1}^{\infty}G_{q^a}$.
  Similarly we define $B_{q^a}$, $T_{q^a}$ and $U_{q^a}$.

\subsection{} Let $N=N_G(T)$ be the normalizer of $T$ in $G$. Then $B$ and $N$ form a $BN$-pair of $G$.
Let $R\subset\text{Hom}(T,\bbbf^*)$ be the root system of $G$ and $R^+$ the set of positive roots determined by $B$. For $\alpha\in R^+$, let $U_\alpha$ be the corresponding root subgroup of $U$.

For any simple root $\alpha$ in $R$, let $s_\alpha$ be the corresponding simple reflection in the Weyl group $W=N/T$. For $w\in W$, $U$ has two subgroups $U_w$ and $U'_w$ such that $U=U'_wU_w$ and $wU'_xw^{-1}\subseteq U$. If $w=s_\alpha$ for some simple root $\alpha$, then $U_w=U_\alpha$ and we simply write $U'_\alpha$ for $U'_w$, which equals $\prod_{\beta\in R^+-\{\alpha\}}U_\beta$. In general,
let $w=s_{\alpha_i}\cdots s_{\alpha_2}s_{\alpha_1}$ be a reduced expression of $w$. Set $\beta_{j}=s_{\alpha_1}s_{\alpha_2}\cdots s_{\alpha_{j-1}}(\alpha_j)$ for $j=1,...,i$. Then

\def\st{\stackrel}
\def\sc{\scriptstyle}

\medskip

\noindent (a) $U_w=U_{\beta_i}\cdots U_{\beta_2}U_{\beta_1}$ and $U'_w=\displaystyle{\prod_{\st{\sc \beta\in R^+}
{w(\beta)\in R^+}}U_\beta}$.

\medskip

\noindent (b) If $\alpha$ and $\beta$ are positive roots  and $w(\alpha)=\beta$, then $n_wU_\alpha n_w^{-1}=U_\beta$, where $n_w$ is a representative of $w$ in $N$.

\medskip

Now assume that $w_0=s_{\alpha_r}\cdots s_{\alpha_2}s_{\alpha_1}$ is a reduced expression of the longest element of $W$. Set $\beta_{j}=s_{\alpha_1}s_{\alpha_2}\cdots s_{\alpha_{j-1}}(\alpha_j)$ for $j=1,...,r$. Then

\medskip

\noindent (c) For any $1\le i\le j\le r$, $U_{\beta_j}\cdots U_{\beta_{i+1}}U_{\beta_i}$ is a subgroup of $U$ and
$U_{\beta_j}\cdots U_{\beta_{i+1}}U_{\beta_i}=U_{\beta_i}U_{\beta_j}\cdots U_{\beta_{i+1}}.$

\medskip

\def\va{\varepsilon_\alpha}

The roots subgroups $U_\alpha,\ \alpha\in R^+$, are also defined over $\bbf$. For each positive root, we fix an isomorphism $\varepsilon_\alpha:\bbbf\to U_\alpha$ such that $t\va(c)t^{-1}=\va(\alpha(t)c)$. Set $U_{\alpha,q^a}=\va(\bbfa)$.

\section{Infinite dimensional Steinberg modules}

In this section the main result (Theorem 2.2) of this paper is proved, which says that certain infinite dimensional Steinberg modules are irreducible.

  \subsection{} Let $k$ be a field. For any one dimensional representation $\theta$ of $T$ over $k$, let $k_\theta$ be the corresponding $kT$-module.
 We define the $k G$-module $M(\theta)=k G\otimes_{k B}k_\theta$. When $\theta$ is trivial representation of $T$ over $k$,  we write $M(tr)$ for $M(\theta)$ and let $1_{tr}$ be a nonzero element in $k_\theta$. We shall also write $x1_{tr}$ instead of $x\otimes 1_{tr}$ for $x\in kG$.

  For $w\in W=N/T$, the element $w1_{tr}$ is defined to be  $n_w1_{tr}$, where $n_w$ is a representative in $N$ of $w$. This is well defined since $T$ acts on $k_\theta$ trivially.
  Let $\eta=\sum_{w\in W} (-1)^{l(w)}w1_{tr}\in M(tr),$ where $l:W\to \mathbb{N}$ is the length function of $W$. Then  $kU\eta$ is a submodule of $M(tr)$ and is called a Steinberg module of $G$, denoted by St, see [X, Prop. 2.3].  Xi proved that St is irreducible if $k$ is the field of complex numbers or $k=\bbbf$ (see [X, Theorem 3.2]).  His argument in fact works for proving that St is irreducible whenever char$k=0$ or char$k=$char$\bbf$. The main result of this paper is the following.

 \subsection{Theorem.}  Assume that $k$ is a field of positive characteristic and char$k\ne$ char$\mathbb{F}_q$. Then the Steinberg module St is irreducible.

Combining Xi's result we have the following result.

 \subsection{Corollary.} The Steinberg module St of $kG$ is irreducible for any field $k$.

 \subsection{} We need some preparation to prove the theorem.

  Let $s_{\alpha_r}s_{\alpha_{r-1}}\cdots s_{\alpha_1}$ be a reduced expression of the longest element $w_0$ of $W$. Set $\beta_i=s_{\alpha_1}\cdots s_{\alpha_{i-1}}(\alpha_i)$. Then $R^+$ consists of $\beta_r,\beta_{r-1},...,\beta_1$, and $U=U_{\beta_r}U_{\beta_{r-1}}\cdots U_{\beta_1}$. Let $n_i$ be a representative in $N=N_G(T)$ of $s_{\alpha_i}$ and set $n=n_rn_{r-1}\cdots n_1$. Note that the elements $z\eta,\ z\in U$, form a basis of St.

 \subsection{Lemma.}  Let $u\in U_{q^a}$. If $u$ is not the neutral element $e$ of $U$, then the sum  of all coefficients of $nu\eta$ in terms the basis $z\eta,\ z\in U$, is 0.

 Proof. Let $\alpha_i$ and $\beta_i$ be as in subsection 2.4. Then $u=u_{r}u_{r-1}\cdots u_1,$ $ u_m\in U_{\beta_m}$. Assume that $u_1=u_2=\cdots u_{i-1}=e$ but $u_i\ne e$, where $e$ is the neutral element of $G$. We use induction on $i$ to prove the lemma. Note that $n_i\eta=-\eta$ for $i=1,2,...,r$.

 Assume that $i=r$. Then $nu\eta=(-1)^{r-1}n_ru'_r\eta$, where $u'_r=n_{r-1}\cdots n_1u_rn_1^{-1}\cdots n_{r-1}^{-1}\in U_{\alpha_r}$. According to the proof of [S, Lemma 1] (see also proof of [X, Proposition 2.3]), there exists $x_r\in U_{\alpha_r}$ such that $n_ru'_r\eta=(x_r-1)\eta$. So the lemma is true in this case.

 Now assume that the lemma is true for $r,r-1,...,i+1$, we show that it is also true for $i$.  In this case we have $nu=(-1)^{i-1}n_r\cdots n_i u'_{r}u'_{r-1}\cdots u'_i\eta$, where $u'_j=n_{i-1}\cdots n_1u_jn_1^{-1}\cdots n_{i-1}^{-1}\in n_{i-1}\cdots n_1U_{\beta_j}n_1^{-1}\cdots n_{i-1}^{-1}= U_{\gamma_j}$, where $\gamma_j=s_{\alpha_i}\cdots s_{\alpha_{j-1}}(\alpha_j)$, $j=i,i+1,...,r$. Then $u'_i\ne e$. Note that $\gamma_i=\alpha_i$. According to the proof of [S, Lemma 1], there exists $x_i\in U_{\gamma_i}$ such that $n_iu'_i\eta=(x_i-1)\eta$.

 If $u'_r\cdots u'_{i+1}=e$, we are done. Now assume that $u'=u'_r\cdots u'_{i+1}\ne e$. Since both $M_i=U_{\gamma_r}U_{\gamma_{r-1}}\cdots U_{\gamma_{i}}$ and $M_{i+1}=U_{\gamma_r}U_{\gamma_{r-1}}\cdots U_{\gamma_{i+1}}$ are subgroups of $U$ and $M_{i+1}x_i= x_iM_{i+1}$,we see that $u'x_i=x_iu''$ for some $u''\in M_{i+1}$ and $u''\ne e$. Thus
 $$(-1)^{i-1}nu\eta=n_r\cdots n_{i+1}u'(x_i-1)\eta=n_r\cdots n_{i+1}x_iu''\eta-n_r\cdots n_{i+1}u'\eta.$$

 By induction hypotheses, we know that the sum of the coefficients of $n_r\cdots n_{i+1}u''\eta$ and the sum of the coefficients of $n_r\cdots n_{i+1}u'\eta$ are 0. Since $n_r\cdots n_{i+1}x_in_{i+1}^{-1}\cdots n_r^{-1}\in U$, we see that the sum of the coefficients of $nu\eta=(-1)^{i-1}n_r\cdots n_{i+1}u'(x_i-1)\eta$ is 0.

 The lemma is proved.

 \subsection{Lemma.}  Let $V$ be a nonzero submodule of St. Then  there exists an integer $a$ such that $\sum_{x\in U_{q^a}}x\eta$ is in $V$.

 Proof. Let $v$ be a nonzero element in $V$. Then $v\in kU_{q^a}\eta$ for some integer $a$. Let $v=\sum_{y\in U_{q^a}}a_yy\eta.$ We may assume that $a_e\ne 0$. Otherwise choose $y\in U_{q^a}$ such that $a_y$ is nonzero and replace $v$ by  $y^{-1}v$.

 By  Lemma 2.5  we see that the sum $A$ of all the coefficients of $nv$ in terms of the basis $z\eta,$ $z\in U$,  is $(-1)^{l(w_0)}a_e\ne 0$.  Thus  $\sum_{\in U_{q^a}}xnv=A\sum_{x\in U_{q^a}} x\eta$. The lemma is proved.

 \subsection{}  Now we can prove the theorem. We show that St=$kGv$ for any nonzero element $v$ in St. Let $V=kGv$.  Let $\alpha_i$ and $\beta_i$ be as in subsection 2.4. For any positive integer $b$, set $X_{i,q^b}=U_{\beta_r,q^b}U_{\beta_{r-1},q^b}\cdots U_{\beta_i,q^b}$. Then $X_{i,q^b}$ is a subgroup of $U$ and $X_{i,q^b}=X_{i+1,q^b}U_{\beta_i,q^b}$. Clearly $X_{i,q^b}$ is a subgroup $X_{i,q^{b'}}$ if $\mathbb{F}_{q^b}$ is a subfield of $\mathbb{F}_{q^{b'}}.$

 We use induction on $i$ to show that there exists positive integer $b_i$ such that the element $\sum_{x\in X_{i,q^{b_i}}}x\eta$ is in $V$. For $i=1$, this is true by Lemma 2.6. Now assume that $\sum_{x\in X_{i,q^{b_i}}}x\eta$ is in $V$, we show that $\sum_{x\in X_{i+1,q^{b_{i+1}}}}x\eta$ is in $V$ for some $b_{i+1}$.

 Let $c_1,...,c_{q^{b_i}+1}$ be a complete set of representatives of all cosets of $\mathbb{F}^*_{q^{b_i}}$ in  $\mathbb{F}^*_{q^{2b_i}}$. Choose $t_1,...,t_{q^{b_i}+1}\in T$ such that $\beta_i(t_j)=c_j$ for $j=1,...,q^{b_i}+1$. Note that $t^{-1}\eta=\eta$ for any $t\in T$. Thus
 $$\sum_{j=1}^{q^{b_i}+1}t_j\sum_{x\in U_{\beta_i,q^{b_i}}}x\eta=q^{b_i}\eta+\sum_{x\in U_{\beta_i,q^{2b_{i}}}}x\eta.$$

 Since $X_{i,q^{b_i}}=X_{i+1,q^{b_i}}U_{\beta_i,q^{b_i}}$ and  $\sum_{x\in X_{i,q^{b_i}}}x\eta$ is in $V$, we see
 \begin{equation*}
 \begin{split}
 \xi=& \sum_{j=1}^{q^{b_i}+1}t_j\sum_{x\in X_{i,q^{b_i}}}x\eta\\
 =&\sum_{j=1}^{q^{b_i}+1}t_j\sum_{y\in X_{i+1,q^{b_i}}}y\sum_{x\in U_{\beta_i,q^{b_i}}}x\eta\\
=&\sum_{y\in X_{i+1,q^{b_i}}}\sum_{j=1}^{q^{b_i}+1}t_jyt_j^{-1}(t_j\sum_{x\in U_{\beta_i,q^{b_i}}}x\eta)\in V.
\end{split}
\end{equation*}

 Choose $b_{i+1}$ such that all $\alpha_m(t_j)$ ($r\ge m\ge i$) are contained in $\mathbb{F}_{q^{b_{i+1}}}$. Then $\mathbb{F}_{q^{b_{i+1}}}$  contains  $\mathbb{F}_{q^{2b_{i}}}$. Thus $t_jyt_j^{-1}$ is in $X_{i+1,q^{b_i+1}}$ for any $y\in X_{i+1,q^{b_i}}$.  Let $Z\in kG$ be the sum of all elements in $X_{i+1,q^{b_i+1}}$. Then we have

 \noindent (1) $\displaystyle {Z\xi=q^{(r-i)b_i}Z\sum_{j=1}^{q^{b_i}+1}t_j\sum_{x\in U_{\beta_i,q^{b_i}}}x\eta=q^{(r-i)b_i}Z(q^{b_i}\eta+\sum_{x\in U_{\beta_i,q^{2b_{i}}}}x\eta)\in V.}$

   Since $\sum_{x\in X_{i,q^{b_i}}}x\eta$ is in $V$, we have $\displaystyle \sum_{x\in X_{i,q^{2b_{i}}}}x\eta\in V.$ Thus

 \noindent(2) $\displaystyle Z\sum_{x\in X_{i,q^{2b_{i}}}}x\eta=q^{2(r-i)b_i}Z\sum_{x\in U_{\beta_i,q^{2b_i}}}x\eta\in V.$

  Since $q\ne 0$ in $k$,  combining (1) and (2) we see that $Z\eta\in V$, i.e., $\sum_{x\in X_{i+1,q^{b_{i+1}}}}x\eta$ is in $V$.

    Note that $X_{r,q^{b_r}}=U_{\beta_r,q^{b_r}}$. Now we have
 $$\sum_{x\in U_{r,q^{b_r}}}x\eta\in V\quad\text{and}\quad\sum_{x\in U_{r,q^{2b_r}}}x\eta\in V.$$
 The above arguments show that
 $$q^{b_r}\eta+\sum_{x\in U_{r,q^{2b_r}}}x\eta\in V.$$
 Therefore $q^{b_r}\eta$ is in $V$. So $V$ contains $kG\eta$=St, hence $V=$St. The theorem is proved.

 \subsection{Remark.} Let St${}_a=kG_{q^a}\eta$. Then St${}_a$ is the Steinberg module of $kG_{q^a}$, which is not irreducible in general. As an example, say, $G=SL_2(\bbbf)$ and $q$ is odd, char$k=2$. Since $q^a+1$ is always divisible by 2=char$k$, St${}_a$ is not irreducible for any positive number $a$ (see [S, Theorems  3]). However, by Theorem 2.2, St is irreducible $kG$-module.

 \section{Non-quasi-finite irreducibility of certain Steinberg modules}

 In this section we show that for certain special linear groups  the Steinberg modules of the groups are not quasi-finite with respect to some natural quasi-finite sequences of the groups, see Proposition 3.2.

 \subsection{} By definition, a group
$G$ is quasi-finite if  $G$ has a sequence $G_1,\ G_2,\ ,...,\ G_n,\
... $ of finite subgroups
 such that $G$ is the union of all $G_i$ and for any positive integers $i,j$ there exists integer $r$
 such that $G_i$ and $G_j$ are contained in $G_r$.  The sequence $G_1,G_2,G_3,...$ is called a quasi-finite sequence of $G$. An irreducible module (or representation) $M$ of
$G$ is {\it quasi-finite} (with respect to the quasi-finite sequence
$G_1,G_2,G_3,...)$ if it has a sequence of subspaces $M_1$, $ M_2$,
$ M_3, $ ... of $M$ such that (1) each $M_i$ is an irreducible
$G_i$-submodule of $M$, (2) if $G_i$ is a subgroup of $G_j$, then
$M_i$ is a subspace of $M_j$, (3) $M$ is the union of all $M_i$. The
sequence $M_1$, $ M_2$, $ M_3, $ ... will be called a quasi-finite
sequence of $M$. See [X, 1.8]

The following question was raised in [4, 1.8]: is every irreducible G-module quasi-finite (with respect to a certain quasi-finite sequence
of G).

The main result of this section is the following result.

\subsection{Proposition.} Let $G=SL_n(\bbbf)$ and $k$ a field  of positive characteristic. Assume that char$k$ divides $(1+q^{a})(1+q^{a}+q^{2a})\cdots(1+q^{a}+\cdots+q^{(n-1)a})$ for all positive integers $a$. If a quasi-finite sequence of $G$  is a subsequence of $SL_n(\bbf),$ $SL_n(\mathbb{F}_{q^2}),$ $ SL_n(\mathbb{F}_{q^3}),$ $ SL_n(\mathbb{F}_{q^4}),...$,  then the Steinberg module St of $kG$ is not quasi-finite with respect to the quasi-finite sequence.

Proof. Let $G_1,\ G_2,\ ,...,\ G_n,\
... $ be a quasi-finite sequence of $G$. Assume that the quasi-finite sequence is a subsequence of $SL_n(\bbf),$ $SL_n(\mathbb{F}_{q^2}),$ $ SL_n(\mathbb{F}_{q^3}),$ $ SL_n(\mathbb{F}_{q^4}),...$. If St is quasi-finite with respect to this quasi-finite sequence, then there exists a sequence of subspaces $M_1$, $ M_2$,
$ M_3, $ ... of St such that (1) each $M_i$ is an irreducible
$G_i$-submodule of $M$, (2) if $G_i$ is a subgroup of $G_j$, then
$M_i$ is a subspace of $M_j$, (3) $M$ is the union of all $M_i$.

Choose a nonzero element $v\in M_1$. By the proof of Theorem 2.2,  there exists $x\in kG$ such that $xv=\eta$. Since $G$ is the union of all $H_a=SL_n(\mathbb{F}_{q^a})$, there exists positive integer $i$ such that $x\in kH_i$. Since $G$ is the union of all $G_a$, there exists $j$ such that $H_i$ is included in $G_j$. Note that $G_j=H_{j'}$ for some positive integer $j'$. Choose integer $d$ such that $G_d$ includes both $G_1$ and $G_j$. Then $M_d$ includes $M_1$ and $M_j$. Moreover, $x\in kG_d$, so that $xv=\eta$ is in $M_d$ and $M_d$ includes $kG_d\eta$. Since $G_d=H_{d'}$ for some $d'$ and $kH_{d'}\eta$ is not irreducible, $M_d$ is not irreducible. This contradicts the assumption that $M_d$ is irreducible $G_d$-module.  The proposition is proved.

\subsection{Remark}  Assume that $n$ is power of a prime $p'$ and $k$ has characteristic $p'$. If $n,q$ are coprime, then $n$ divides $(1+q^{a})(1+q^{a}+q^{2a})\cdots(1+q^{a}+\cdots+q^{(n-1)a})$  for all positive integers $a$. To see this, let $A_m= 1+q^{a}+\cdots+q^{(m-1)a}$, $m=2,...,n$. Since $n,q$ are coprime, if $A_m\equiv 1$(mod $n$), then $m>2$ and $A_{m-1}$ is divisible by $n$. If $A_m\not\equiv 1$(mod $n$) for $m=2,...,n$, then either some $A_m$ is divisible by $n$ or $A_m\equiv A_l$(mod $n$) for some $n\ge m>l\ge 2$. Since  $n,q$ are coprime, we have $m\ge l+2$. Then $A_{m-l}$ is divisible by $n$. By Proposition 2.2, in this case, the Steinberg module St of $kSL_n(\bbbf)$ is not quasi-finite for a quasi-finite sequence of $G$ whenever it is a subsequence of $SL_n(\bbf),$ $SL_n(\mathbb{F}_{q^2}),$ $ SL_n(\mathbb{F}_{q^3}),$ $ SL_n(\mathbb{F}_{q^4}),...$.

However, it is not clear that whether  St is quasi-finite with other quasi-finite sequences of $SL_n(\bbbf)$.

For other reductive groups, one discusses similarly.

\subsection{} Assume that $G$ is quasi-finite and has sequence of normal subgroups  $\{1\}=G_0\subset G_1\subset \cdots\subset
G_n=G$ such that all
$G_i/G_{i-1}$ are abelian. Xi asks whether any irreducible $\mathbb{C}G$-module is isomorphic to the induced module of a one dimensional
module of a subgroup of $G$ (see [X, 1.12]). The question has a negative answer ever for finite groups, for instance, the two-dimensional irreducible complex representation of $SL_2(\mathbb F_3)$ is an counterexample. Perhaps for the question the condition of all
$G_i/G_{i-1}$ being abelian should be further strengthened to all
$G_i/G_{i-1}$ being cyclic.

\bigskip

{\bf Acknowledgement.} I am very grateful to Professor Nanhua Xi for  guidance and  great helps in writing the paper. The work was done during my visit to the Academy of Mathematics and Systems Science, Chinese Academy of Sciences as a visiting student for the academic year 2014-2015. I thank the  Academy of Mathematics and Systems Science for hospitality.

\bigskip


\end{document}